# The Number of Different Effective Partitions and a Specific Goldbach Partition of Any Given Even Number Greater Than *6*

Song Linggen
(Fudan Univ., Shanghai (200433), P. R. China)

**Keywordsv**
Effective primes/products/integers/partitions of even number *2N* greater than *6*;
The number of different/identical effective partitions of *2N*;
Relatively prime theorem;
Two-part method;
A specific Goldbach partition of any given even number greater than *6*.

**Abstract**

Other than any odd prime whose factor is contained by the given even number $2N_m$ greater than *6*, total *m (m ≥ 2)* odd primes within open interval *(1, $2N_m$ - 1)* were defined as effective primes $p_{(i, m)}$ of $2N_m$, where *i = 1, 2, …, m - 1, m*, and $p_{(1, m)}$ and $p_{(m, m)}$ is the minimum and maximum one among $p_{(i, m)}$, respectively.

Let *Eq.s* $\alpha_{(i, m)}$ & *Eq.s* $\beta_{(i, m)}$ represent *2m* simultaneous equations $p_{(i, m)} + \alpha_{(i, m)} = 2N_m$ and $p_{(1, m)} p_{(i, m)} + \beta_{(i, m)} = 2N_m$, respectively.

By using two-part method starting from *m = 2*, it was verified that whatever *m* is, the number of different effective partition(s) of $2N_m$ among $2N_m$ simultaneous *Eq.s* $\alpha_{(i, m)}$ & *Eq.s* $\beta_{(i, m)}$ is always less than *m*, and consequently any given even number $2N_m$ greater than *6* must be divided by $p_{(g, m)}$ into a Goldbach partition of $2N_m$, where $p_{(g, m)}$ represents the greatest one among all the effective primes of $2N_m$ whose factors are contained by $\alpha_{(i, m)}$.

A specific Goldbach partition of any given even number greater than *6* can be found definitely.

**Text**
Enlightened by Ben Green and Terence Tao[1], this manuscript is to establish a way to definitely find a specific Goldbach partition of any given even number greater than *6*.

## §1, Basic concepts
### 1-1, Goldbach partition of even number *2N* greater than *6*
When a given even number *2N* greater than *6* can be partitioned to two odd primes, such partition is called a Goldbach partition of *2N*.

### 1-2, Effective primes/products/integers/partitions of the given even number $2N_m$
Other than any odd prime whose factor is contained by the given even number $2N_m$ greater than 6, total *m* odd primes within open interval *(1, $2N_m$ - 1)* were defined as effective primes $p_{(i, m)}$ of $2N_m$ ordered by magnitudes as $1 < p_{(1, m)} < p_{(2, m)} < ... < p_{(m-1, m)} < p_{(m, m)} < 2N_m - 1$.

Within open interval *(1, $2N_m$ - 1)*, products which contain and only contain the effective prime factor(s) of $2N_m$ were defined as effective products of $2N_m$.

Effective primes and effective products of $2N_m$ were commonly named as effective integers of $2N_m$. It is obvious that every effective integer of $2N_m$ is within open interval *(1, $2N_m$ - 1)*.

The partition of $2N_m$ which consists of two effective integers of $2N_m$ was defined as an effective partition of $2N_m$.

### 1-3, The number of effective primes of any given even number $2N_m$ greater than *6* is always greater than *1*
Any given even number $2N_m$ greater than 6 is corresponding to a definite number *m* of effective primes of $2N_m$, and *m* is always no less than 2 because effective integers $(N_m - k)$ and $(N_m + k)$ of $2N_m$ are relatively prime to each other, where $k = 1$ or $2$ when $N_m$ is even or odd, respectively.

## §2, Theorem *1*
Theorem *1* states that:

Let

$p_{(i, m)} + \alpha_{(i, m)} = 2N_m$            Eq.s $\alpha_{(i, m)}$
$p_{(1, m)} p_{(i, m)} + \beta_{(i, m)} = 2N_m$       Eq.s $\beta_{(i, m)}$

where, $p_{(i, m)}$ represent total *m* effective primes of the given even number $2N_m$ greater than 6;
    $i = 1, 2, ..., m - 1, m$, respectively;
    $m \geqslant 2$;
    $p_{(1, m)}$ and $p_{(m, m)}$ is the minimum and maximum one among $p_{(i, m)}$, respectively;

Each of $\alpha_{(i, m)}$ is an effective integer of $2N_m$, respectively;

$\beta_{(i, m)}$ is an effective integer of $2N_m$ only when $p_{(1, m)}p_{(i, m)} < 2N_m - 1$, respectively.

Whatever $m$ is, the number of different effective partition(s) of $2N_m$ among $2m$ simultaneous equations $Eq.s\ \alpha_{(i, m)}$ & $Eq.s\ \beta_{(i, m)}$ is always less than $m$.

### § 3, Some necessary preparative concepts for the verification of Theorem *1*
### 3-1, Properties of $\alpha_{(i, m)}$ & $\beta_{(i, m)}$

Because of relatively prime theorem and the definitions of effective primes/integers of $2N_m$, there are *5* useful properties of $\alpha_{(i, m)}$ & $\beta_{(i, m)}$ for the verification of Theorem *1*.

*Prop. 1*

Each of $\alpha_{(i, m)}$ is an effective integer without factor $p_{(i, m)}$ within open interval $(1, 2N_m - 1)$;

$2N_m - 1 > \alpha_{(1, m)} > \alpha_{(2, m)} > \cdots > \alpha_{(m-1, m)} > \alpha_{(m, m)} \geq p_{(1, m)}$;

$\beta_{(i, m)}$ is an effective integer without factors $p_{(i, m)}$ and $p_{(1, m)}$ only when $p_{(1, m)}p_{(i, m)} < 2N_m - 1$;

$2N_m - 1 > \beta_{(1, m)} > \beta_{(2, m)} > \cdots > \beta_{(m-1, m)} > \beta_{(m, m)}$;

$p_{(1, m)}p_{(m, m)} > p_{(1, m)}p_{(m-1, m)} > \cdots > p_{(1, m)}p_{(2, m)} > p_{(1, m)}p_{(1, m)}$.

*Prop. 2*

There is no any common divisor between $\alpha_{(i, m)}$ and $\beta_{(i, m)}$, respectively.

*Prop. 3*

At most one of $\alpha_{(i, m)}$ and $\beta_{(i, m)}$ may contain the maximum effective prime factor $p_{(m, m)}$ of $2N_m$, respectively.

*Prop. 4*

$\beta_{(m, m)}$ is not an effective integer of $2N_m$ as long as there is a Goldbach partition of $2N_m$.

*Prop. 5*

There are at most $(m - k)/2$ successive $\beta_{(j, m)}$ equals an effective prime of $2N_m$, respectively, where, $j = 1, 2, \ldots, (m - k)/2$, respectively; $k = 0$ or $1$ when $m$ is even or odd, respectively.

### 3-2, The number of different partition(s) of $2N_m$

There are totally $2m$ simultaneous equations $Eq.s\ \alpha_{(i, m)}$ & $Eq.s\ \beta_{(i, m)}$.

On the following *3* conditions only, the number $D_m$ of different effective partitions of $2N_m$ among these $2m$ simultaneous equations is reduced from $2m$:

*Condition 1*

When $\alpha_{(x, m)} = p_{(y, m)}$,

$p_{(x, m)} + p_{(y, m)} = 2N_m$          Eq. $\alpha_{(x, m)}$
$p_{(y, m)} + p_{(x, m)} = 2N_m$          Eq. $\alpha_{(y, m)}$

Eq. $\alpha_{(x, m)}$ and Eq. $\alpha_{(y, m)}$ are identical, and $D_m$ is one reduced from $2m$;

*Condition 2*
When $\beta_{(x, m)} = p_{(y, m)}$,

$p_{(1, m)}p_{(x, m)} + p_{(y, m)} = 2N_m$          Eq. $\beta_{(x,, m)}$
$p_{(y, m)} + p_{(1, m)}p_{(x, m)} = 2N_m$          Eq. $\alpha_{(y, m)}$

Eq. $\beta_{(x, m)}$ and Eq. $\alpha_{(y, m)}$ are identical, and $D_m$ is one reduced from $2m$;

*Condition 3*
When $p_{(1, m)}p_{(z, m)} \geq 2N_m - 1$,

because $1 \geq \beta_{(z, m)} > \beta_{(z+1, m)} > \cdots > \beta_{(m-1, m)} > \beta_{(m, m)}$, that is Eq. $\beta_{(z, m)}$, Eq. $\beta_{(z+1, m)}$, ..., Eq. $\beta_{(m-1, m)}$ and Eq. $\beta_{(m, m)}$ are not effective partitions of $2N_m$, and $D_m$ is $(m - z + 1)$ reduced from $2m$.

The number of different partition(s) of $2N_m$ among total $2m$ simultaneous equations Eq.s $\alpha_{(i, m)}$ & Eq.s $\beta_{(i, m)}$ could be expressed as the following *Eq.(3, 1)*:

$D_m = 2m - a_m/2 - b_m - c_m$          Eq.(**3, 1**)

where, $D_m$ represents the number of different effective partition(s) of $2N_m$ among total $2m$ simultaneous equations Eq.s $\alpha_{(i, m)}$ & Eq.s $\beta_{(i, m)}$;
    $a_m$ among $m$ $\alpha_{(i, m)}$ equals an effective prime of $2N_m$, respectively;
    $a_m$ is a non-odd integer;
    $b_m$ among $m$ $\beta_{(i, m)}$ equals an effective prime of $2N_m$, respectively;
    $b_m \geq 0$;
    $c_m$ among $m$ $\beta_{(i, m)}$ is not an effective integer of $2N_m$, respectively;
    $c_m \geq 0$;

$$m \geqslant b_m + c_m.$$

### 3-3 $D_m$, $D'_m$ and $d_m$

On other hand, $D_m$ could also be expressed as the following *Eq.(3, 2)*:

$$D_m = D'_m + d_m \qquad Eq.(3, 2)$$

where, $D_m$ represents the number of different effective partition(s) of $2N_m$ among total $2m$ simultaneous equations *Eq.s* $\alpha_{(i, m)}$ & *Eq.s* $\beta_{(i, m)}$, where $i = 1, 2, \ldots, m - 1, m$, respectively;

$D'_m$ represents the number of different effective partition(s) of $2N_m$ in which the maximum effective prime factor $p_{(m, m)}$ of $2N_m$ not take part among total $2(m-1)$ simultaneous equations *Eq.s* $\alpha_{(i, m)}$ & *Eq.s* $\beta_{(i, m)}$, where $i = 1, 2, \ldots, m - 2, m - 1$, respectively;

$d_m$ represents the number of different effective partition(s) of $2N_m$ in which the maximum effective prime factor $p_{(m, m)}$ of $2N_m$ takes part among total $2m$ equations *Eq.s* $\alpha_{(i, m)}$ & *Eq.s* $\beta_{(i, m)}$, where $i = 1, 2, \ldots, m - 1, m$, respectively;

$m \geqslant 2$.

### 3-4 Relationship between $D_{m-1}$ and $D'_m$

According to the similarity between the definitions of $D_{m-1}$ and $D'_m$, the following *In-eq.(3, 3)* holds:

$$D'_m \leqslant w \text{ when } D_{m-1} \leqslant w \qquad In\text{-}eq.(3, 3)$$

where, $w$ represents an integer no less than $1$.

$m \geqslant 2$.

### §4, Verification of Theorem *1*

Based on these preparative concepts mentioned above and by using two-part method starting from $m = 2$, Theorem *1* was verified as shown below.

### 4-1, Theorem *1* holds when $m = 2$

There are total two effective primes, $p_{(1, 2)}$ and $p_{(2, 2)}$, of the given even number $2N_2$ when $m = 2$.

*4* simultaneous equations *Eq.s* $\alpha_{(i, 2)}$ & *Eq.s* $\beta_{(i, 2)}$ were listed below:

$p_{(1, 2)} + a_{(1, 2)} = 2N_2$    Eq. $a_{(1, 2)}$         $p_{(1, 2)}p_{(1, 2)} + \beta_{(1, 2)} = 2N_2$    Eq. $\beta_{(1, 2)}$

$p_{(2, 2)} + a_{(2, 2)} = 2N_2$    Eq. $a_{(2, 2)}$         $p_{(1, 2)}p_{(2, 2)} + \beta_{(2, 2)} = 2N_2$    Eq. $\beta_{(2, 2)}$

According to *Prop.s 1* and *2*, $a_{(1, 2)}$ and $a_{(2, 2)}$ has to equal an exact power of prime factor $p_{(2, 2)}$ and $p_{(1, 2)}$, respectively, and both $\beta_{(1, 2)}$ and $\beta_{(2, 2)}$ contains neither factor $p_{(1, 2)}$ nor factor $p_{(2, 2)}$.

According to the definition of effective integers of $2N_2$, both $\beta_{(1, 2)}$ and $\beta_{(2, 2)}$ are not effective integers of $2N_2$, and the following *In-eq.(4, 1)* holds:

$1 \geqslant \beta_{(1, 2)} > \beta_{(2, 2)}$        *In-eq.(4, 1)*

Because $p_{(1, 2)}p_{(2, 2)} > p_{(1, 2)}p_{(1, 2)} \geqslant 2N_2 - 1$ when *In-eq.(4, 1)* holds, $a_{(1, 2)}$ and $a_{(2, 2)}$ has to equal prime $p_{(2, 2)}$ and $p_{(1, 2)}$, respectively.

Therefore, *Eq.s* $a_{(i, 2)}$ & *Eq.s* $\beta_{(i, 2)}$ could be revised as the following:

$p_{(1, 2)} + p_{(2, 2)} = 2N_2$    Eq. $a_{(1, 2)}$         $p_{(1, 2)}p_{(1, 2)} + B_{(1, 2)} = 2N_2$    Eq. $\beta_{(1, 2)}$

$p_{(2, 2)} + p_{(1, 2)} = 2N_2$    Eq. $a_{(2, 2)}$         $p_{(1, 2)}p_{(2, 2)} + B_{(2, 2)} = 2N_2$    Eq. $\beta_{(2, 2)}$

where, both $B_{(1, 2)}$ and $B_{(2, 2)}$ are not effective integers of $2N_2$.

Observing the revised *Eq.s* $a_{(i, 2)}$ & *Eq.s* $\beta_{(i, 2)}$ gave that $a_2 = 2$, $b_2 = 0$ and $c_2 = 2$. According to *Eq.(3, 1)*, $D_2 = 1$ when $m = 2$, and consequently Theorem *1* holds when $m = 2$.

On other hand, observing the revised *Eq.s* $a_{(i, 2)}$ & *Eq.s* $\beta_{(i, 2)}$ gave that $D'_2 = 0$ and $d_2 = 1$. According to *Eq.(3, 2)*, $D_2 = 1$ when $m = 2$, and consequently Theorem *1* holds when $m = 2$.

Therefore, Theorem *1* holds when $m = 2$.

### 4-2 Theorem *1* holds when $m = M +1$ if Theorem *1* holds when $m = M$   ($M \geqslant 2$)

*Eq.s* $a_{(i, M+1)}$ & *Eq.s* $\beta_{(i, M+1)}$ when $m = M +1$ were listed below:

$p_{(i, M+1)} + a_{(i, M+1)} = 2N_{M+1}$    *Eq.s* $a_{(i, M+1)}$        $p_{(1, M+1)}p_{(i, M+1)} + \beta_{(i, M+1)} = 2N_{M+1}$    *Eq.s* $\beta_{(i, M+1)}$

where, $p_{(i, M+1)}$ represents total $M+1$ effective primes of the given even number $2N_{M+1}$;
  $i = 1, 2, \ldots, M, M+1$, respectively;
  $M \geqslant 2$;
  $p_{(1, M+1)}$ and $p_{(M+1, M+1)}$ is the minimum and maximum one among $p_{(i, M+1)}$, respectively;
  Each of $\alpha_{(i, M+1)}$ is an effective integer of $2N_{M+1}$, respectively;
  $\beta_{(i, M+1)}$ is an effective integer of $2N_{M+1}$ only when $p_{(1, M+1)}p_{(i, M+1)} < 2N_{M+1} - 1$.

It was verified that $D_{M+1} \leqslant (M+1) - 1$ if $D_M \leqslant M - 1$ as shown below:

### 4-2-1 $D'_{M+1} \leqslant M - 1$ if Theorem 1 holds when $m = M$

Assuming Theorem 1 holds when $m = M$ implies that the following *In-eq.*(**4, 2**) holds:

$$D_M \leqslant M - 1 \qquad \textit{In-eq.}(\textbf{4, 2})$$

where, $D_M$ represents the number of different effective partitions of $2N_M$ among $2M$ equations
  *Eq.s* $\alpha_{(i, M)}$ & *Eq.s* $\beta_{(i, M)}$;
  $i = 1, 2, \ldots, M - 1, M$, respectively.

According to *In-eq.*(**3, 3**) and *In-eq.*(**4, 2**), the following *In-eq.*(**4, 3**) holds, too:

$$D'_{M+1} \leqslant M - 1 \text{ when } D_M \leqslant M - 1 \qquad \textit{In-eq.}(\textbf{4, 3})$$

where, $D'_{M+1}$ represents the number of different effective partitions of $2N_{M+1}$ in which only factors $p_{(i, M+1)}$ take part among $2M$ simultaneous equations *Eq.s* $\alpha_{(i, M+1)}$ & *Eq.s* $\beta_{(i, M+1)}$;
  $D_M$ represents the number of different effective partitions of $2N_M$ in which only factors $p_{(i, M)}$ take part among $2M$ simultaneous equations *Eq.s* $\alpha_{(i, M)}$ & *Eq.s* $\beta_{(i, M)}$;
  $i = 1, 2, \ldots, M - 1, M$, respectively.

### 4-2-2, $d_{M+1} = 1$ if Theorem 1 holds when $m = M$

The following *In-eq.*(**4, 4**) was derived from *Eq.*(**3, 1**), *Eq.*(**3, 2**) and *In-eq.*(**4, 3**)

$$a_{M+1}/2 + b_{M+1} + c_{M+1} + d_{M+1} \geqslant M+3 \quad \text{when } D_M \leqslant M - 1 \qquad \textit{In-eq.}(\textbf{4, 4})$$

where, $a_{M+1}$ among $M+1$ $\alpha_{(i, M+1)}$ equals an effective prime of $2N_{M+1}$, respectively;
$a_{M+1}$ is a non-odd integer;
$b_{M+1}$ among $M+1$ $\beta_{(i, M+1)}$ equals an effective prime of $2N_{M+1}$, respectively;
$b_{M+1} \geqslant 0$;
$c_{M+1}$ among $M+1$ $\beta_{(i, M+1)}$ is not an effective integer of $2N_{M+1}$, respectively;
$c_{M+1} \geqslant 0$;
$d_{M+1}$ represents the number of different effective partition(s) of $2N_{M+1}$ in which the factor $p_{(M+1, M+1)}$ of $2N_{M+1}$ takes part among total $2(M+1)$ simultaneous equations Eq.s $\alpha_{(i, M+1)}$ & Eq.s $\beta_{(i, M+1)}$;
$M+1 \geqslant b_{M+1} + c_{M+1}$.
$i = 1, 2, …, M, M+1$, respectively.

According to *Prop. 3*, in terms of how $\alpha_{(i, M+1)}$ and $\beta_{(i, M+1)}$ contain factor $p_{(M+1, M+1)}$, there are totally *4* sub-conditions in principle when $m = M+1$.

It was verified that the following *In-eq.(4, 5)* always holds on the each logical sub-conditions:

$d_{M+1} = 1$ when $D_M \leqslant M - 1$      *In-eq.(**4**, **5**)*

### Sub-condition 1, None of $\alpha_{(i, M+1)}$ and $\beta_{(i, M+1)}$ contained factor $p_{(M+1, M+1)}$
$d_{M+1} \leqslant 2$, because factor $p_{(M+1, M+1)}$ only takes part in the following 2 simultaneous equations:

$p_{(M+1, M+1)} + \alpha_{(M+1, M+1)} = 2N_{M+1}$      Eq. $\alpha_{(M+1, M+1)}$
$p_{(1, M+1)}p_{(M+1, M+1)} + \beta_{(M+1, M+1)} = 2N_{M+1}$      Eq. $\beta_{(M+1, M+1)}$

<u>*1, When $p_{(1, M+1)}p_{(M+1, M+1)} \geqslant 2N_{M+1} – 1$*</u>
Eq. $\beta_{(M+1, M+1)}$ is not an effective partition of $2N_{M+1}$, and only Eq. $\alpha_{(M+1, M+1)}$ is an effective partition of $2N_{M+1}$, $d_{M+1}=1$;

<u>*2, If $p_{(1, M+1)}p_{(M+1, M+1)} \leqslant 2N_{M+1} - 1$*</u>
According to *Prop. 4*, there would be no any Goldbach partition of $2N_{M+1}$, $a_{M+1} = 0$;

Each of $\beta_{(i, M+1)}$ would be an effective integer of $2N_{M+1}$, respectively, $c_{M+1} = 0$;

Because $a_{M+1} = 0$, $c_{M+1} = 0$, and $d_{M+1} \leqslant 2$, observing *In-eq.(**4**, **4**)* would gave *In-eq.(**4**, **6**)*:

$b_{M+1} \geq M+1$  when $D_M \leq M - 1$       In-eq.(**4, 6**)

According to *Prop. 5*, this is illogical.

Therefore, because $p_{(1, M+1)}p_{(M+1, M+1)} \geq 2N_{M+1} - 1$ always holds, $d_{M+1} = 1$ on **Sub-condition 1**.

### Sub-condition 2, Only $\alpha_{(a, M+1)}$ contains factor $p_{(M+1, M+1)}$

$d_{M+1} \leq 3$, because factor $p_{(M+1, M+1)}$ only takes part in the following 3 simultaneous equations:

$p_{(a, M+1)} + Q_a p_{(M+1, M+1)}^{f(a)} = 2N_{M+1}$       Eq. $\alpha_{(a, M+1)}$
$p_{(M+1, M+1)} + \alpha_{(M+1, M+1)} = 2N_{M+1}$       Eq. $\alpha_{(M+1, M+1)}$
$p_{(1, M+1)}p_{(M+1, M+1)} + \beta_{(M+1, M+1)} = 2N_{M+1}$       Eq. $\beta_{(M+1, M+1)}$

where, $p_{(a, M+1)}$ is an effective prime of $2N_{M+1}$;
  $a \neq M+1$.
  $Q_a$ is an effective integer of $2N_{M+1}$ or $1$;
  Exponent $f(a)$ is an integer no less than $1$.

**1, When $p_{(1, M+1)}p_{(M+1, M+1)} \geq 2N_{M+1} - 1$**
Eq. $\beta_{(M+1, M+1)}$ is not an effective partition of $2N_{M+1}$;

According to *Prop. 1*, effective integer $\alpha_{(a, M+1)}$ of $2N_{M+1}$ has to equal effective prime $p_{(M+1, M+1)}$ of $2N_{M+1}$, and therefore Eq. $\alpha_{(a, M+1)}$ and Eq. $\alpha_{(M+1, M+1)}$ are identical.

There is only one different effective partition of $2N_{M+1}$ among Eq. $\alpha_{(a, M+1)}$, Eq. $\alpha_{(M+1, M+1)}$ and Eq. $\beta_{(M+1, M+1)}$.

$d_{M+1} = 1$.

**2, If $p_{(1, M+1)}p_{(M+1, M+1)} \leq 2N_{M+1} - 1$**
According to *Prop. 4*, there would be no any Goldbach partition of $2N_{M+1}$, $a_{M+1} = 0$;

Each of $\beta_{(i, M+1)}$ would be an effective integer of $2N_{M+1}$, respectively, $c_{M+1} = 0$;

Because $a_{M+1} = 0$, $c_{M+1} = 0$, and $d_{M+1} \leq 3$, observing In-eq.(**4, 4**) would gave In-eq.(**4, 7**):

$b_{M+1} \geq M$   when $D_M \leq M - 1$       In-eq.(**4, 7**)

According to *Prop. 5*, this is illogical.

Therefore, because $p_{(1, M+1)}p_{(M+1, M+1)} \geq 2N_{M+1} - 1$ always holds, $d_{M+1} = 1$ on **Sub-condition 2**.

**Sub-condition 3, Only $\beta_{(b, M+1)}$ contains factor $p_{(M+1, M+1)}$**

$d_{M+1} \leq 3$, because factor $p_{(M+1, M+1)}$ only takes part in the following 3 simultaneous equations:

$$p_{(1, M+1)}p_{(b, M+1)} + Q_b p_{(M+1, M+1)}^{f(b)} = 2N_{M+1} \quad \text{Eq. } \beta_{(b, M+1)}$$
$$p_{(M+1, M+1)} + \alpha_{(M+1, M+1)} = 2N_{M+1} \quad \text{Eq. } \alpha_{(M+1, M+1)}$$
$$p_{(1, M+1)}p_{(M+1, M+1)} + \beta_{(M+1, M+1)} = 2N_{M+1} \quad \text{Eq. } \beta_{(M+1, M+1)}$$

where, $p_{(b, M+1)}$ is an effective prime of $2N_{M+1}$;
$b \neq M+1$.
$Q_b$ is an effective integer of $2N_{M+1}$ or $1$;
Exponent $f(b)$ is an integer no less than $1$.

<u>1, When $p_{(1, M+1)}p_{(M+1, M+1)} \geq 2N_{M+1} - 1$</u>
Eq. $\beta_{(M+1, M+1)}$ is not an effective partition of $2N_{M+1}$;

Eq. $\beta_{(b, M+1)}$ is not an effective partition of $2N_{M+1}$ unless $Q_b p_{(M+1, M+1)}^{f(b)} = p_{(M+1, M+1)}$;

When $Q_b p_{(M+1, M+1)}^{f(b)} = p_{(M+1, M+1)}$, Eq. $\beta_{(b, M+1)}$ and Eq. $\alpha_{(M+1, M+1)}$ are identical.

There is only one different effective partition of $2N_{M+1}$ among Eq. $\beta_{(b, M+1)}$, Eq. $\alpha_{(M+1, M+1)}$ and Eq. $\beta_{(M+1, M+1)}$.

$d_{M+1} = 1$.

<u>2, If $p_{(1, M+1)}p_{(M+1, M+1)} < 2N_{M+1} - 1$</u>
According to *Prop. 4*, there would be no any Goldbach partition of $2N_{M+1}$, $a_{M+1} = 0$;

Each of $\beta_{(i, M+1)}$ would be an effective integer of $2N_{M+1}$, respectively, $c_{M+1} = 0$;

Because $a_{M+1} = 0$, $c_{M+1} = 0$, and $d_{M+1} \leq 3$, observing *In-eq.(4, 4)* would gave *In-eq.(4, 8)*:

$$b_{M+1} \geq M \quad \text{when } D_M \leq M - 1 \qquad \text{In-eq.}(4, 8)$$

According to *Prop. 5*, this is illogical.

Therefore, because $p_{(1, M+1)} p_{(M+1, M+1)} \geq 2N_{M+1} - 1$ always holds, $d_{M+1} = 1$ on **Sub-condition 3**.

**Sub-condition 4**, Only $\alpha_{(a, M+1)}$ and $\beta_{(b, M+1)}$ contain factor $p_{(M+1, M+1)}$

$d_{M+1} \leq 4$, because factor $p_{(M+1, M+1)}$ only takes part in the following 4 simultaneous equations:

| | |
|---|---|
| $p_{(a, M+1)} + Q_a p_{(M+1, M+1)}^{f(a)} = 2N_{M+1}$ | Eq. $\alpha_{(a, M+1)}$ |
| $p_{(1, M+1)} p_{(b, M+1)} + Q_b p_{(M+1, M+1)}^{f(b)} = 2N_{M+1}$ | Eq. $\beta_{(b, M+1)}$ |
| $p_{(M+1, M+1)} + \alpha_{(M+1, M+1)} = 2N_{M+1}$ | Eq. $\alpha_{(M+1, M+1)}$ |
| $p_{(1, M+1)} p_{(M+1, M+1)} + \beta_{(M+1, M+1)} = 2N_{M+1}$ | Eq. $\beta_{(M+1, M+1)}$ |

where, $p_{(a, M+1)}$ is an effective prime of $2N_{M+1}$;
    $a \neq M+1$.
    $Q_a$ is an effective integer of $2N_{M+1}$ or $1$;
    Exponent $f(a)$ is an integer no less than $1$.
    $p_{(b, M+1)}$ is an effective prime of $2N_{M+1}$;
    $b \neq M+1$.
    $Q_b$ is an effective integer of $2N_{M+1}$ or $1$;
    Exponent $f(b)$ is an integer no less than $1$.

Observing Eq. $\alpha_{(a, M+1)}$ and Eq. $\beta_{(b, M+1)}$ gave that if this sub-condition existed, the following Eq.(**4, 9**) and Eq.(**4, 10**) would have to hold:

| | |
|---|---|
| $Q_a p_{(M+1, M+1)}^{f(a)} = p_{(1, M+1)} p_{(M+1, M+1)}$ | Eq.(**4, 9**) |
| $Q_b p_{(M+1, M+1)}^{f(b)} = p_{(M+1, M+1)}$ | Eq.(**4, 10**) |

The 4 simultaneous equations in which factor $p_{(M+1, M+1)}$ takes part were revised below:

| | |
|---|---|
| $p_{(a, M+1)} + p_{(1, M+1)} p_{(M+1, M+1)} = 2N_{M+1}$ | Eq. $\alpha_{(a, M+1)}$ |
| $p_{(1, M+1)} p_{(b, M+1)} + p_{(M+1, M+1)} = 2N_{M+1}$ | Eq. $\beta_{(b, M+1)}$ |
| $p_{(M+1, M+1)} + p_{(1, M+1)} p_{(b, M+1)} = 2N_{M+1}$ | Eq. $\alpha_{(M+1, M+1)}$ |
| $p_{(1, M+1)} p_{(M+1, M+1)} + p_{(a, M+1)} = 2N_{M+1}$ | Eq. $\beta_{(M+1, M+1)}$ |

Observing these revised equations gave that Eq. $\alpha_{(a, M+1)}$ and Eq. $\beta_{(M+1, M+1)}$ would be identical;

Observing these revised equations gave that *Eq.* $\alpha_{(M+1, M+1)}$ and *Eq.* $\beta_{(b, M+1)}$ would be identical;

This implies that $d_{M+1} \leqslant 2$, and **_Sub-condition 4_** would be similar to **_Sub-condition 1_**.

Because $p_{(1, M+1)}p_{(M+1, M+1)} \geqslant 2N_{M+1} - 1$ always holds if $d_{M+1} \leqslant 2$ as verified on **_Sub-condition 1_**, **_Sub-condition 4_** is illogical.

Summarizing above gave that $d_{M+1} = 1$ holds when $m = M + 1$ if Theorem *1* holds when $m = M$.

Because $d_{M+1} = 1$ if Theorem *1* holds when $m = M$, the following *In-eq.*(*4, 11*) was derived from *Eq.*(*3, 2*) and *In-eq.*(*4, 3*):
:

$D_{M+1} \leqslant (M+1) - 1$ when $D_M \leqslant M - 1$          *In-eq.*(*4, 11*)

Because Theorem *1* holds when $m = 2$ as verified in **_4-1_**, according to *In-eq.*(*4, 11*), whatever $m$ is, Theorem *1* always holds.

### *§5*, Theorem *2* and its verification
#### *5-1*, **Theorem *2***
Theorem *2* states that:

Any given even number $2N_m$ greater than *6* must be divided by its effective prime $p_{(g, m)}$ into a Goldbach partition of $2N_m$, where $p_{(g, m)}$ represents the greatest one among all the effective primes of $2N_m$ whose factors are contained by $\alpha_{(i, m)}$.

#### *5-2*, **Verification of Theorem *2***
According to Theorem *1*, there are at most *(m - 1)* different effective partitions of $2N_m$ among $m$ simultaneous equations *Eq.s* $\alpha_{(i, m)}$ ($i = 1, 2, …, m - 1, m$, respectively). Therefore, there is at least one Goldbach partition of $2N_m$:

$p_{(a, m)} + p_{(a', m)} = 2N_m$          *Eq.*(*5, 1*)

where, $p_{(a, m)}$ and $p_{(a', m)}$ are two effective primes of $2N_m$.

The following *In-eq.(5, 2)*, *In-eq.(5, 3)* and *In-eq.(5, 4)* always holds:

$$p_{(m, m)} \geqslant p_{(g, m)} \geqslant p_{(a, m)} > N_m > p_{(a', m)} \geqslant p_{(1, m)} \quad \text{In-eq.(5, 2)}$$
$$p_{(1, m)} \geqslant 3 \quad \text{In-eq.(5, 3)}$$
$$p_{(1, m)} p_{(g, m)} > 2N_m - 1 \quad \text{In-eq.(5, 4)}$$

where, $p_{(1, m)}$ and $p_{(m, m)}$ is the minimum and maximum effective prime of $2N_m$, respectively;
  $p_{(g, m)}$ is the greatest one among all the effective primes of $2N_m$ whose factors take part in $a_{(i, m)}$;
  $i = 1, 2, \ldots, m-1, m$, respectively.

Let $a_{(g', m)}$ contains factor $p_{(g, m)}$:

$$p_{(g', m)} + Q_{g'} \cdot p_{(g', m)}^{f(g')} = 2N_m \quad \text{Eq. } a_{(g', m)}$$

where, $Q_{g'}$ is an effective integer of $2N_m$ or $1$;
  Exponent $f(g')$ is an integer no less than $1$.

Because of *In-eq.(5, 4)*, there is no any effective product containing factor $p_{(g, m)}$ within open interval $(1, 2N_m - 1)$, and *Eq. $a_{(g', m)}$* can be revised as the following:

$$p_{(g', m)} + p_{(g, m)} = 2N_m \quad \text{Eq. } a_{(g', m)}$$

where, both $p_{(g', m)}$ and $p_{(g, m)}$ are effective primes of $2N_m$.

Therefore, Theorem 2 always holds.

## §6, Conclusion
A specific Goldbach partition of any given even number greater than 6 can be found definitely.

## The End